\documentclass[12pt]{amsart}

\newtheorem{theorem}{Theorem}[section]
\newtheorem{lemma}[theorem]{Lemma}
\newtheorem{proposition}[theorem]{Proposition}

\theoremstyle{definition}

\newtheorem{remark}[theorem]{Remark}

\usepackage{amscd,amssymb}

\begin{document}

\title[Gr\"obner bases of ideals invariant under endomorphisms]
{Gr\"obner bases of ideals \\
 invariant under endomorphisms}

\author[Vesselin Drensky and Roberto La Scala]
{Vesselin Drensky and Roberto La Scala}
\address{Institute of Mathematics and Informatics,
Bulgarian Academy of Sciences,
Acad. G. Bonchev Str., Block 8, 1113 Sofia, Bulgaria}
\email{drensky@math.bas.bg}

\address{Dipartimento di Matematica,
Universit\`a di Bari, Via E. Orabona 4,
70125 Bari, Italia}
\email{lascala@dm.uniba.it}

\thanks
{The research of the first author was partially supported by Grant
MI-1503/2005 of the Bulgarian National Science Fund.}
\thanks{The research of the second was partially supported
by Universit\`a di Bari.}


\begin{abstract}
We introduce the notion of Gr\"obner $S$-basis of an ideal of the free
associative algebra $K\langle X\rangle$ over a field $K$
invariant under the action of a semigroup $S$ of endomorphisms
of the algebra. We calculate the Gr\"obner $S$-bases of the ideal corresponding to
the universal enveloping algebra of the free nilpotent of class 2 Lie algebra
and of the T-ideal generated by the polynomial identity $[x,y,z]=0$,
with respect to suitable semigroups $S$.
In the latter case, if $\vert X\vert>2$,
the ordinary Gr\"obner basis is infinite and our Gr\"obner $S$-basis is finite.
We obtain also explicit minimal Gr\"obner bases of these ideals.
\end{abstract}
\maketitle

\section{Introduction}
Let $K$ be a field of any characteristic
and let $X = \{x_1,x_2,\ldots\}$ be a finite or countable
set with more than one element.
Let $K\langle X\rangle$ be the free unitary associative
$K$-algebra generated by $X$. Its elements are
polynomials in the noncommuting variables $x_i$.

In this paper we study some two-sided ideals of $K\langle X\rangle$
from computational point of view. We immediately face the problem
that, even when the set $X$ is finite, many important ideals
of $K\langle X\rangle$ are not finitely generated. On the other hand,
quite often these ideals have additional
structure and ``uniformly looking'' generating sets.

For example, let $L=L(X)$ be the free Lie algebra freely generated
by $X$ and canonically embedded into $K\langle X\rangle$.
It is known that the free nilpotent of class $c$ Lie algebra
$L/[\underbrace{L,\ldots,L}_{c+1\text{ times}}]$,
usually denoted in the theory of varieties of Lie algebras
$F({\mathfrak N}_c)$,
has a set of defining relations
consisting of all (left normed) commutators
\[
u_j=[[\ldots [x_{j_1},x_{j_2}],\ldots,x_{j_c}],x_{j_{c+1}}]=0.
\]
Hence, by the Poincar\'e-Birkhoff-Witt theorem,
its universal enveloping algebra $U(F({\mathfrak N}_c))$
is a homomorphic image of $K\langle X\rangle$ modulo the ideal
$I$ generated by all $u_j$. We may define the ideal $I$
as the minimal ideal of $K\langle X\rangle$ which contains
the commutator
\[
[x_1,x_2,\ldots,x_c,x_{c+1}]
=[[\ldots [x_1,x_2],\ldots,x_c],x_{c+1}]
\]
and is invariant, or stable, under all endomorphisms sending $X$ to $X$.

Other examples are the T-ideals of $K\langle X\rangle$.
These ideals are invariant under all endomorphisms of
$K\langle X\rangle$ and coincide with the ideals of polynomial
identities of suitable PI-algebras. If
\[
U=\{u_j(x_1,\ldots,x_{n_j})\mid j\in J\}\subset K\langle X\rangle
\]
is any set, then the T-ideal generated by $U$ is generated
as a usual ideal by all $u_j(f_1,\ldots,f_{n_j})$, when
the $n_j$ polynomials $f_1,\ldots,f_{n_j}$ run on $K\langle X\rangle$.
For infinite $X$ nontrivial T-ideals cannot be finitely generated as ideals.
If the set $X$ is finite, then a theorem of Markov \cite{Ma}
describes the few cases when a T-ideal is finitely generated as an ideal.
This happens if and only if the T-ideal contains
for some $c$ the Engel polynomial
\[
[x_2,\underbrace{x_1,\ldots,x_1}_{c\text{ times}}].
\]
One of the classical problems in PI-theory is the Specht problem
\cite{Sp} which states whether any T-ideal is finitely generated
as a T-ideal. The celebrated structure theory of T-ideals developed by Kemer,
see his book \cite{K2} for the account, allowed him \cite{K1} in 1987
to give a positive solution to the Specht problem over
a field of characteristic 0.
In the case of positive characteristic there are several counterexamples.
The first of them were given by Belov \cite{Be}, Grishin \cite{G},
and Shchigolev \cite{Sh}.
A good source for the state of the art of the Specht problem,
as well as an improved exposition of the theory of Kemer, can be found
in the recent book by Kanel-Belov and Rowen \cite{KBR}.

When the set of variables $X$ is finite,
the knowledge of a generating set of an ideal $I$ of
the polynomial algebra $K[X]$ is not always sufficient for
concrete calculations with the elements of $I$
and in the factor algebra $K[X]/I$. A similar phenomenon
appears for the ideals of $K\langle X\rangle$, even if the
ideal has a finite generating set.
In commutative algebra the problem is solved with the technique of
Gr\"obner bases. This is a powerful tool for computing with
commutative algebras, in algebraic geometry, and in invariant theory,
see e.g. the books by Adams and Loustaunau \cite{AL}, Kreuzer and Robbiano
\cite{KRo}, Sturmfels \cite{St}. One may introduce Gr\"obner bases for ideals
not only in the polynomial algebra $K[X]$, but also in the free associative algebra
$K\langle X\rangle$ and in the free Lie algebra $L(X)$. Shirshov \cite{S} proved
his Composition Lemma dealing with Lie polynomials
via associative words. His algorithms for free Lie algebras work also for free
associative algebras. In the noncommutative case one often
calls the corresponding bases Gr\"obner-Shirshov bases instead
of Gr\"obner bases.
In the last three decades the number of the applications of the
noncommutative Gr\"obner bases  increases, see e.g. the seminal papers
by Bokut \cite{B} and Bergman \cite{Bg}, the surveys by Mora \cite{Mo} and Ufnarovski \cite{U},
as well as the relatively recent surveys by
Bokut, Fong, Ke, and Kolesnikov \cite{BFKK}, and
Bokut and Kolesnikov {BK1, BK2}.
Nevertheless, there are very few examples of ideals of the free algebras
with explicitly known Gr\"obner bases. Also, it is a well known fact that
many algorithmic problems are not solvable
for associative algebras, with the word problem among them, and in the general
case there is no algorithm to construct a Gr\"obner basis
of an ideal of a free associative algebra. Recently, Gr\"obner-Shirshov
bases have been introduced also for other
free objects, see e.g. Bokut, Fong, and Shiao \cite{BFS}.

In the present paper we consider ideals $I$ of the free algebra
$K\langle X\rangle$ which are invariant under the action of a subsemigroup
$S$ of the endomorphism semigroup of $K\langle X\rangle$. We introduce
the notion of Gr\"obner $S$-basis of $I$. This is a subset $B$ of $I$
with the property that $S(B)$ is a Gr\"obner basis of $I$
in the usual sense, with respect to some term-ordering of the monomials
in $K\langle X\rangle$.

We handle completely two cases of Gr\"obner $S$-bases. The first is
the universal enveloping algebra $U(F({\mathfrak N}_2))$
of the free nilpotent of class 2
Lie algebra $F({\mathfrak N}_2)=L/[L,L,L]$. The semigroup $S$
consists of all endomorphisms which send $X$ to $X$ and preserve
the ordering on $X$. The corresponding ideal $I$ of $K\langle X\rangle$
is generated by the commutators $[x_i,x_j,x_k]$. We give a concrete
finite Gr\"obner $S$-basis of $I$. It consists of commutators of length 3
and one more commutator of degree 4.

As Lalonde and Ram \cite{LR},
and Bokut and Malcolmson \cite{BM}, mentioned,
if $H$ is an ideal of the free Lie algebra $L(X)$ and $B$ is its
Gr\"obner-Shishov basis with respect to a certain ordering on
a suitable basis of the vector space $L(X)$, then $B$ is also
a Gr\"obner basis of the ideal $I$ of $K\langle X\rangle$ generated
by $H$. Of course, the factorization modulo this ideal gives the
universal enveloping algebra $U(L/H)$. This result easily implies that
the algebra $U(F({\mathfrak N}_2))$ does have a Gr\"obner basis
consisting of polynomials of degree 3 and 4 only. We want to mention
that our approach is direct and does not use the facts from \cite{LR}
and \cite{BM}. Instead, we use easy combinatorics of words
and the explicit $K$-basis of $U(F({\mathfrak N}_2))$.

The second example treats another algebra of importance for the theory of PI-algebras
and with applications to superalgebras. This is the relatively free
algebra $F(\text{\rm var}E)$
of the variety of associative algebras generated by the Grassmann (or exterior)
algebra $E$ over an infinite field of characteristic different from 2.
This algebra can be considered as the generic Grassmann algebra.
The structure of $F(\text{\rm var}E)$, $\text{\rm char}K=0$,
was described by Krakowski and Regev \cite{KR}, see also
the paper by Di Vincenzo \cite{DV} or the book \cite{D} by one of the authors.
It is known that the polynomial identities of $E$ are consequences of
the commutator identity $[x_1,x_2,x_3]=0$. The defining relations of $F(\text{\rm var}E)$
in characteristic 0 were described by Latyshev \cite{L} and
consist of the polynomials
\[
[x_i,x_j,x_k]=0,\quad
[x_i,x_j][x_k,x_l]+[x_i,x_k][x_j,x_l]=0,
\]
where $x_i,x_j,x_k,x_l$ are replaced by all possible elements of $X$.
It is well known that the same polynomials form a set of defining relations
of $K\langle X\rangle/([x_1,x_2,x_3])^T$ over any field of
characteristic different from 2, where $([x_1,x_2,x_3])^T$ is the T-ideal of $K\langle X\rangle$
generated by $[x_1,x_2,x_3]$.
Bokut and Makar-Limanov \cite{BML} showed that, when $\vert X\vert>2$,
the ideal $([x_1,x_2,x_3])^T$
has no finite Gr\"obner basis. On the other hand, they introduced an extra set of
generators of the algebra $F(\text{var}E)$, $y_{ij}=[x_i,x_j]$, which are in its centre, and established that
the corresponding Gr\"obner basis is finite when $X$ is finite.
In the present paper we show that, although the Gr\"obner basis of the
T-ideal $([x_1,x_2,x_3])^T$ is infinite for $m>2$, it is uniformly looking.
We present explicitly a finite set of polynomials $G$ and a subsemigroup $S$ of the endomorphism semigroup
of $K\langle X\rangle$ such that $G$ is a Gr\"obner $S$-basis of the ideal.
We also correct some inaccuracies in the paper by Bokut and Makar-Limanov  \cite{BML}.
Again, our approach is based on combinatorics of words and the explicit basis of $F(\text{var}E)$.

\section{$S$-ideals, $S$-bases, Gr\"obner $S$-bases}

Denote by $\text{End}(K\langle X\rangle)$ the semigroup of all endomorphisms
of the $K$-algebra $K\langle X\rangle$. Let $S\subset \text{End}(K\langle X\rangle)$ be
a subsemigroup which includes the identity endomorphism.
If $I$ is a two-sided ideal of $K\langle X\rangle$ we say that $I$
is an {\em $S$-invariant ideal} or simply an {\em $S$-ideal} if it is
invariant under all the endomorphisms of $S$, i.e.
\[
\varphi(I) \subset I \quad \text{for all }\varphi\in S.
\]

To construct an $S$-ideal it is sufficient to take any subset $B\subset K\langle X\rangle$
and form the two-sided ideal $I$ generated by $S(B)$. In this case, we say
that $B$ is an {\em $S$-basis} of $I$.

A natural problem is to establish if, for different choices of the semigroup $S$,
all the $S$-ideals have finite $S$-bases. For example, the positive solution
by Kemer \cite{K1} of the Specht problem in characteristic 0 can be restated
that for $S = \text{End}(K\langle X\rangle)$ every $S$-invariant ideal
is finitely $S$-generated.

We fix now on $K\langle X\rangle$ a {\em term-ordering} $<$, i.e.
a linear order on the set $\langle X\rangle$ of words, or monomials,
which is a multiplicatively compatible well-ordering.
This means that:

(i) For every two different monomials $u,v$ we have either
$u<v$ or $v<u$;

(ii) Every subset of $\langle X\rangle$ has a minimal element;

(iii) If $u<v$ in $\langle X\rangle$, then $wu<wv$ and $uw<vw$
for every $w\in \langle X\rangle$.

If $f\in K\langle X\rangle$ is a nonzero polynomial,
we denote by $\text{\rm lt}(f)$
the greatest monomial of $f$. We recall that a
{\em Gr\"obner basis} of an ideal $I$ of $K\langle X\rangle$
is a subset $G\subset I$ (not necessarily finite) which satisfies
the following property: for any nonzero $f\in I$ there exists
a $g_i\in G$ such that $\text{lt}(g_i)$ is a subword of $\text{lt}(f)$.
By induction on the term-ordering, it is easy to prove that we can
write any $f\in I$ as
\[
f = \sum f_i g_i h_i,
\]
where $g_i\in G$ (possibly $g_i = g_j$ for $i\neq j$) and we have
$f_i,h_i\in K\langle X\rangle$, only a finite
number of them different from zero, such that for all $i$
\[
\text{lt}(f) \geq \text{lt}(f_i)\text{lt}(g_i)\text{lt}(h_i).
\]
We have hence that $G$ is also a generating set of $I$ as a two-sided
ideal of $K\langle X\rangle$. For any subset
$G\subset K\langle X\rangle$ it is useful to define
$\text{init}(G)$ as the two-sided ideal
generated by the set of monomials
$\{\,\text{lt}(g_i)\mid g_i\in G\,\}$. We say
that $\text{init}(G)$ is the {\em initial ideal} generated by $G$.
Then, we have clearly that a subset $G\subset I$ is a Gr\"obner basis
of $I$ if and only if $\text{init}(G) = \text{init}(I)$.
In other words, the set of monomials of $K\langle X\rangle$
\[
\{\,w\in\langle X\rangle \mid \text{there exists } g_i\in G\ \text{such that}\
\text{lt}(g_i)\ \text{is a subword of}\ w\,\}
\]
is a $K$-basis of the subspace $\text{init}(I)\subset K\langle X\rangle$.
Then the set
\[
N=\langle X\rangle\backslash\text{lt}(\text{init}(G))
\]
of {\em normal words} with respect to $G$ is
a $K$-basis of the factor algebra $K\langle X\rangle/I$.
A Gr\"obner basis $G$ of an ideal $I$ is called
{\em reduced} if every $g_i\in G$ is a linear combination of
normal words with respect to $G\backslash \{g_i\}$.
Moreover, we call $G$ {\em minimal} if for any $g_i\in G$
we have that $G\backslash \{g_i\}$ is not a Gr\"obner basis of $I$,
i.e. $\text{lt}(g_i)$ is a normal word with respect to
$G\backslash \{g_i\}$. For more details about the theory of
noncommutative Gr\"obner bases we refer to \cite{Mo, U}.

Now let $S$ be a semigroup of endomorphisms of $K\langle X\rangle$ and
let $G$ be a subset of the $S$-ideal $I$. We say that $G$ is a
{\em Gr\"obner $S$-basis} of $I$ if $S(G)$ is a Gr\"obner basis of $I$
as a two-sided ideal of $K\langle X\rangle$, i.e. $\text{init}(I)$
is equal to the initial ideal generated by $S(G)$.

\section{Universal enveloping algebras of free nilpotent algebras}

We keep $K$, $X$, and $K\langle X\rangle$ as in the previous section.
We introduce the standard {\em deg-lex} ordering on $\langle X\rangle$.
We compare the monomials first by total degree
and then lexicographically, reading them from left to right
and assuming that $x_1<x_2<\cdots$. We consider
$K\langle X\rangle$ also as a multigraded vector space, counting
in the monomials the number of enterings of each variable.
If $f,g\in K\langle X\rangle$, the {\em commutator} of $f,g$ is simply
the polynomial
\[
[f,g] = f g - g f.
\]
We refer to the book by Bahturin \cite{Ba} as a background
on Lie algebras and their polynomial identities. Here we summarize
the basic facts we need.
The Lie subalgebra of $K\langle X\rangle$ generated by $X$ with respect
to the commutator operation is isomorphic to the free Lie algebra
freely generated by $X$. We denote this algebra by $L=L(X)$.
Every Lie algebra generated by a countable (or finite) set
is isomorphic to $L/H$ for some ideal $H$ of the Lie algebra $L$.
Then the Poincar\'e-Birkhoff-Witt theorem gives that the universal
enveloping algebra $U(L/H)$ is isomorphic to
$K\langle X\rangle/I$, where $I=K\langle X\rangle HK\langle X\rangle$
is the ideal of $K\langle X\rangle$ generated by $H$. If
$f_1,f_2,\ldots$ is a basis of the $K$-vector space $L/H$, then
$U(L/H)$ has a $K$-basis consisting of all ``monomials''
$f_1^{a_1}\cdots f_p^{a_p}$.

The algebra $L$ has several important bases consisting of commutators.
They are built on the
following principle. One fixes an ordered set of associative
{\em Lyndon-Shirshov monomials}
defined in terms of some special combinatorial properties.
Then, for each monomial, one arranges the Lie brackets in a certain
recursive way, and obtains the basis of $L$.
The elements of the basis are either elements of $X$ or commutators
$[[u],[v]]$, where $[u],[v]$ are also elements of the basis.
The bases under consideration
allow to introduce an analogue of Gr\"obner bases for the ideals $H$ of $L$,
called {\em Gr\"obner-Shirshov} bases, see the original paper by Shirshov \cite{S}.
The subset $G$ of $H$ is a Gr\"obner-Shirshov basis of $H$, if for every
nonzero $f\in H$ with leading commutator $[u]$ there exists a $g\in G$
with leading commutator $[v]$ such that the associative word $v$ obtained by
deleting the Lie brackets in $[v]$ is a subword of the associative word $u$.
As we discussed in the introduction, every Gr\"obner-Shirshov
basis of the ideal $H$ of $L$ is a Gr\"obner basis of the ideal
$I$ generated in $K\langle X\rangle$ by $H$.

We denote by $F({\mathfrak N}_c)$ the free nilpotent of class $c$ Lie algebra.
It is isomorphic to the factor algebra of $L$ modulo the
$(c+1)$-st member
\[
\gamma_{c+1}(L)=[\underbrace{L,\ldots,L}_{c+1 \text{ times}}]
\]
of the lower central series of $L$. It is well known that $\gamma_{c+1}(L)$
is spanned by all commutators of length $\geq c+1$
and can be generated as an ideal by commutators of length $c+1$.
The following easy statement is well known and we omit the proof.

\begin{proposition}\label{result of Bokut and Malcolmson}
There exists a Gr\"obner basis with respect to the deg-lex ordering
of the ideal of $K\langle X\rangle$ generated by $\gamma_{c+1}(L)$
consisting only of commutators of length $c+1,c+2,\ldots,2c$.
\end{proposition}

We apply Proposition \ref{result of Bokut and Malcolmson} to
the ideal $I$ of $K\langle X\rangle$ generated by $\gamma_3(L)$.

\begin{proposition}\label{Groebner basis for gamma3}
The polynomials
\begin{equation}\label{commutators of length 3 in two variables}
f'_{ij}=[[x_i,x_j],x_j],\quad f''_{ij}=[x_i,[x_i,x_j]],\quad i>j,
\end{equation}
\begin{equation}\label{commutators of length 3 in three variables}
g'_{ijk}=[x_i,[x_j,x_k]],\quad g''_{ikj}=[[x_i,x_k],x_j],\quad i>j>k,
\end{equation}
\begin{equation}\label{commutators of length 4}
h_{ijk} = [[x_i,x_j],[x_i,x_k]],\quad i>j>k,
\end{equation}
form a Gr\"obner basis with respect to the deg-lex ordering
of the ideal $I = K\langle X\rangle \gamma_3(L) K\langle X\rangle$.
\end{proposition}

\begin{proof}
We consider the set $B$ of all associative Lyndon-Shirshov words
$u$ defined with the property
that $u$ is bigger than all its cyclic rearrangements. The brackets on $u$
are arranged as follows.
One finds the longest right Lyndon-Shirshov subword $v$ of $u$.
Then $u=wv$ for some word $w$. It turns out that $w$ is also a Lyndon-Shirshov word.
One considers the nonassociative Lyndon-Shirshov words $[w],[v]$ corresponding to $w$ and $v$.
Then one defines $[u]=[[w],[v]]$.

By Proposition \ref{result of Bokut and Malcolmson}
we need all associative Lyndon-Shirshov words of length 3 and 4. They are
\begin{equation}\label{Lyndon-Shirshov words of degree 3 in 2 variables}
x_ix_jx_j,\quad x_ix_ix_j,\quad i>j,
\end{equation}
\begin{equation}\label{Lyndon-Shirshov words of degree 3 in 3 variables}
x_ix_jx_k,\quad x_ix_kx_j,\quad i>j>k,
\end{equation}
\begin{equation}\label{Lyndon-Shirshov words of degree 4 in 2 variables}
x_ix_jx_jx_j,\quad x_ix_ix_jx_j,\quad x_ix_ix_ix_j,\quad i>j,
\end{equation}
\begin{equation}\label{Lyndon-Shirshov words of degree 4 in 3 variables}
x_ix_ix_jx_k,\quad i>j,k,\quad x_ix_jx_ix_k,\quad i>j>k,
\end{equation}
\begin{equation}\label{Lyndon-Shirshov words of degree 4 in 4 variables}
x_ix_jx_kx_l,\quad i>j,k,l.
\end{equation}
The arrangement of the brackets in the cases
(\ref{Lyndon-Shirshov words of degree 3 in 2 variables})
and (\ref{Lyndon-Shirshov words of degree 3 in 3 variables}) is, respectively,
\[
[[x_i,x_j],x_j],\quad [x_i,[x_i,x_j]],\quad
[x_i,[x_j,x_k]],\quad [[x_i,x_k],x_j],
\]
and this gives the elements $f'_{ij}$ and $f''_{ij}$, $i>j$,
in (\ref{commutators of length 3 in two variables}) and
$g'_{ijk}$ and $g''_{ijk}$, $i>j>k$, in
(\ref{commutators of length 3 in three variables}).
Similarly, we obtain $h_{ijk}$, $i>j>k$, in
(\ref{commutators of length 4}) from $x_ix_jx_ix_k$
in (\ref{Lyndon-Shirshov words of degree 4 in 3 variables}).
We do not need the commutators built on the words
from (\ref{Lyndon-Shirshov words of degree 4 in 2 variables}),
(\ref{Lyndon-Shirshov words of degree 4 in 4 variables}),
and the words $x_ix_ix_jx_k$ from
(\ref{Lyndon-Shirshov words of degree 4 in 3 variables})
for the Gr\"obner basis of the ideal $I$ generated
by $\gamma_3(L)$ because they contain a subword of the form
$u_iu_iu_j$ or $u_iu_ju_k$ with $i>j,k$.
Hence the commutators  (\ref{commutators of length 3 in two variables}),
(\ref{commutators of length 3 in three variables}), and
(\ref{commutators of length 4}) give a Gr\"obner basis of $I$.
\end{proof}

Now we state Proposition \ref{Groebner basis for gamma3}
in terms of Gr\"obner $S$-bases.

\begin{theorem}\label{Groebner S-basis for nilpotent Lie algebras}
Let $X$ be an infinite set and let $S$ be the semigroup
consisting of all endomorphisms of $K\langle X\rangle$ which send $X$ to $X$ and preserve
the ordering on $X$. Then the set of polynomials
\begin{equation}\label{Groebner S-basis for gamma3 A}
[[x_2,x_1],x_1],\quad [x_2,[x_2,x_1]],
\end{equation}
\begin{equation}\label{Groebner S-basis for gamma3 B}
[x_3,[x_2,x_1]],\quad [[x_3,x_1],x_2],\quad [[x_3,x_2],[x_3,x_1]]
\end{equation}
is a Gr\"obner $S$-basis of the ideal of $K\langle X\rangle$
generated by $\gamma_3(L)$.
\end{theorem}

\begin{proof}
Let $\varphi_1$ be an endomorphism from $S$ such that $\varphi_1(x_1)=x_j$
and $\varphi_1(x_2)=x_i$, $i>j$. Applying $\varphi_1$ to $[[x_2,x_1],x_1]$ and $[x_2,[x_2,x_1]]$
we obtain the elements (\ref{commutators of length 3 in two variables}).
Similarly, if $i>j>k$, we start with $\varphi_2\in S$ satisfying
 $\varphi_2(x_1)=x_k$,  $\varphi_2(x_2)=x_j$, $\varphi_2(x_3)=x_i$. Applying it on
$[x_3,[x_2,x_1]]$, $[[x_3,x_1],x_2]$, and $[[x_3,x_2],[x_3,x_1]]$,
we obtain (\ref{commutators of length 3 in three variables})
and (\ref{commutators of length 4}).
In this way, acting by $S$ on the elements from
(\ref{Groebner S-basis for gamma3 A}) and (\ref{Groebner S-basis for gamma3 B}),
we obtain the Gr\"obner basis of the ideal generated by $\gamma_3(L)$.
\end{proof}

\begin{remark}\label{X is finite}
(i) It is easy to see that applying the semigroup $S$ from Theorem
\ref{Groebner S-basis for nilpotent Lie algebras} to the Gr\"obner $S$-basis
(\ref{Groebner S-basis for gamma3 A}), (\ref{Groebner S-basis for gamma3 B}),
we obtain a minimal Gr\"obner basis which is not reduced. The polynomial
\[
[x_3,[x_2,x_1]]=x_3x_2x_1-x_3x_1x_2-x_2x_1x_3+x_1x_2x_3
\]
contains as a summand the monomial $x_3x_1x_2$ which can be reduced using
$[[x_3,x_1],x_2]$. The commutator $[[x_3,x_2],[x_3,x_1]]$ also needs to be
reduced. These reductions can be done by easy calculations.

(ii) The restriction that $X$ is countable is not essential.
Theorem \ref{Groebner S-basis for nilpotent Lie algebras} can be restated
for any infinite well-ordered set $X$.

(iii)  When the set $X$ is finite, the semigroup $S$ from
Theorem \ref{Groebner S-basis for nilpotent Lie algebras} consists of the identity endomorphism only.
We may replace it with the semigroup generated by the
endomorphisms $\varphi_1$, $\varphi_2$ of $K\langle X\rangle$
with the property $\varphi_1(X),\varphi_2(X)\subseteq X$, $\varphi_1(x_1)<\varphi_1(x_2)$,
$\varphi_2(x_1)<\varphi_2(x_2)<\varphi_2(x_3)$.
\end{remark}

We shall give another direct combinatorial description of the Gr\"obner basis
of the ideal of $K\langle X \rangle$ generated by $\gamma_3(L)$ which
we shall use later for the Gr\"obner basis of the T-ideal $([x_1,x_2,x_3])^T$.

\begin{lemma}\label{basis of universal enveloping algebra}
The polynomials
\begin{equation}\label{basis of length one and two}
x_{i_1}\cdots x_{i_l}[x_{j_1},x_{k_1}]\cdots[x_{j_m},x_{k_m}] ,
\end{equation}
where  $i_1\leq\cdots\leq i_l$, $j_s>k_s$, $s=1,\ldots,m$,
and $(j_1,k_1)\leq\cdots\leq (j_m,k_m)$
with respect to the lexicographic ordering,
form a $K$-basis of the universal enveloping algebra $U(F({\mathfrak N}_2))$.
\end{lemma}

\begin{proof}
The Poincar\'e-Birkhoff-Witt theorem
gives that, if $g_1,g_2,\ldots$ is an ordered $K$-basis of a Lie algebra, then its universal
enveloping algebra has a $K$-basis consisting of all $g_1^{a_1}\cdots g_n^{a_n}$.
This immediately completes the proof:
the free nilpotent of class 2 Lie algebra $F({\mathfrak N}_2)$ is spanned by
all commutators of length 1 and 2, i.e. by the elements $x_i$ and $[x_i,x_j]$,
and the anticommutativity allows to assume that $i>j$ in $[x_i,x_j]$.
\end{proof}

\begin{lemma}\label{normal words for gamma3}
The set of normal words $N(G)$ with respect to the set $G$ of the commutators
(\ref{commutators of length 3 in two variables}),
(\ref{commutators of length 3 in three variables}), and
(\ref{commutators of length 4}) consists of all monomials
$w = x_{i_1} \cdots x_{i_n}$ such that
\begin{itemize}
\item[(i)] The inequality $i_k > i_{k+1}$ implies that $i_k \leq i_{k+2}$ and if, additionally
$k>1$, then $i_{k-1}<i_k$;
\item[(ii)] If $i_k = i_{k+2} > i_{k+1},i_{k+3}$, then
$i_{k+1}\leq i_{k+3}$.
\end{itemize}
\end{lemma}

\begin{proof}
The leading monomials of the elements of $G$ are of three types:
\begin{itemize}
\item[(a)] $x_ix_jx_j$ and $x_ix_ix_j$, where $i\geq j$;
\item[(b)] $x_ix_jx_k$, where $i>j,k$;
\item[(c)] $x_i x_j x_i x_k$, where $i > j > k$.
\end{itemize}
If the word $w = x_{i_1} \cdots x_{i_n}$ does not satisfy
the condition (i), then $i_k > i_{k+1}$ for some $k$, but $i_{k-1} \geq i_k$
or $i_k > i_{k+2}$. In this case at least one of the subwords $x_{i_{k-1}} x_{i_k} x_{i_{k+1}}$
and $x_{i_k} x_{i_{k+1}} x_{i_{k+2}}$ is of type (a) or (b). Suppose now that
$w$ does not satisfy (ii), i.e. $i_k = i_{k+2} > i_{k+1},i_{k+3}$
and $i_{k+1} > i_{k+3}$. Then, the subword
$x_{i_k} x_{i_{k+1}} x_{i_{k+2}} x_{i_{k+3}}$ is of type (c).
Moreover, the above arguments can be clearly reversed.
\end{proof}

Now we give an explicit bijection between the basis of
$U(F({\mathfrak N}_2))$ from Lemma
\ref{basis of universal enveloping algebra} and the set of
normal words from Lemma \ref{normal words for gamma3}.

\begin{proposition}\label{correspondence}
There is a one-to-one correspondence between the set $B$
of the products (\ref{basis of length one and two}) and the set
$N(G)$ from Lemma \ref{normal words for gamma3} which preserves the multigrading.
\end{proposition}

\begin{proof}
Although the statement of the proposition is almost obvious by
the construction of the considered algebras and ideals, we give a formal proof.
We consider the set of sequences of indices that parameterize the polynomials in $B$,
say:
\[
\overline{B} = \{ (i_1,\ldots,i_l,(j_1,k_1),\ldots,(j_m,k_m)) \}.
\]
We consider also the set of sequences of indices that occur
in the words of $N=N(G)$:
\[
\overline{N} = \{ (i_1,\ldots,i_n)\mid i_k\ \text{ satisfies (i),(ii)} \}
\]
We define recursively a map $\psi$ from $\overline{B}$ into
the set of sequences of integers. If
$u = (i_1,\ldots,i_l,(j_1,k_1),\ldots,(j_m,k_m))$ then we find the first index $i_{p+1}$ with the
property  $j_1 \leq i_{p+1}$ (hence $i_p<j_1$ if $p\geq 1$) and define
\[
\psi(u) = (i_1,\ldots,i_p,j_1,k_1,\psi(v))
\]
where $v = (i_{p+1},\ldots,i_l,(j_2,k_2),\ldots,(j_m,k_m))$. We shall prove that
the image of $\psi$ is contained in $\overline{N}$. Since $i_1\leq\cdots\leq i_l$
and by the definition of $\psi$ we have that $\psi(u)$ satisfies
the condition (i). Moreover, owing to the lexicographic ordering of the
pairs $(j_1,k_1),\ldots,(j_m,k_m)$ it is clear that also (ii) is verified.
For example, if
\[
u = (1,2,2,2,3,4,5,6,({\bf 2},{\bf 1}),({\bf 2},{\bf 1}),({\bf 3},{\bf 1}),
({\bf 3},{\bf 2}),({\bf 5},{\bf 2}),({\bf 5},{\bf 3}),({\bf 6},{\bf 4})),
\]
(we have typesetted the pairs $(j,k)$ in bold) then
\begin{equation}\label{psi of u}
\psi(u)=(1,{\bf 2},{\bf 1},{\bf 2},{\bf 1},2,2,2,{\bf 3},{\bf 1},{\bf 3},{\bf 2},
3,4,{\bf 5},{\bf 2},{\bf 5},{\bf 3},5,{\bf 6},{\bf 4},6).
\end{equation}

We define now two maps $\vartheta_1,\vartheta_2$ from $\overline{N}$ respectively into
the set of integer sequences and the set of sequences of pairs of integers.
If $u = (i_1,\ldots,i_n)$ then:
\[
\vartheta_1(u) = (i_1,\ldots,i_{k-1},\vartheta_1(v))\ \text{and}\
\vartheta_2(u) = ((i_k,i_{k+1}),\vartheta_2(v)),
\]
where $i_1 \leq \cdots \leq i_k > i_{k+1}$ and
$v = (i_{k+2},\ldots,i_n)$. Define now the map
$\vartheta:v \mapsto (\vartheta_1(v),\vartheta_2(v))$. We claim that the image of $\vartheta$
is contained in $\overline{B}$. In fact, by the condition (i) we have that
$\vartheta_1(u)$ is an increasing sequence of indices. Moreover, from the
definition of $\vartheta_2$ and the condition (ii) it follows that $\vartheta_2(u)$
is a sequence of pairs $(j,k)$ with $j > k$, which is increasing with
respect to the lexicographic ordering. In the above example, if $v=\psi(u)$,
then $\vartheta(v)=u$.

Finally, it is easy to check that the maps $\psi$ and $\vartheta$ induce bijections
between $B$ and $N$ which preserve the multigrading and are
inverse of each other.
\end{proof}

\begin{remark}
It is more convenient, compare with the example in (\ref{psi of u}), to write the
normal words $N(G)$ from Lemma \ref{normal words for gamma3} in the form
\begin{equation}\label{easy form for normal words}
x_1^{a_1}(x_2x_1)^{b_{21}}x_2^{a_2}(x_3x_1)^{b_{31}}(x_3x_2)^{b_{32}}x_3^{a_3}\cdots
\prod_{p=1}^{m-1}(x_mx_p)^{b_{mp}}x_m^{a_m},
\end{equation}
where $a_i,b_{ij}\geq 0$. For example, in (\ref{psi of u}) we have the word
\[
x_1(x_2x_1)^2x_2^3(x_3x_1)(x_3x_2)x_3x_4(x_5x_2)(x_5x_3)x_5(x_6x_4)x_6.
\]
\end{remark}

Let $I$ be a multigraded ideal of $K\langle X\rangle$ and let $B$ be a multigraded basis
of $R=K\langle X\rangle/I$. If $G$ is a subset of $I$ and $N(G)$ is the set of normal
words with respect to $G$, then in each multihomogeneous component of $B$ and $N(G)$,
the number of elements from $B$ is not greater than the number of elements from $N(G)$. If
the number of these elements coincides for each multihomogeneous component, we have
that $G$ is a Gr\"obner basis for $I$. Hence
Proposition \ref{correspondence} implies immediately
Proposition \ref{Groebner basis for gamma3}
and Theorem \ref{Groebner S-basis for nilpotent Lie algebras}.

\section{The polynomial identities of the Grassmann algebra}

In this section we assume that the base field $K$ is of characteristic different from 2.
We consider the T-ideal $T=([x_1,x_2,x_3])^T$ of $K\langle X\rangle$
generated by the commutator $[x_1,x_2,x_3]$. We shall summarize the necessary facts,
including also some proofs to make the exposition self-contained.
The following proposition is well known, see the paper by Latyshev \cite{L}
or the book \cite{D} for the case of characteristic 0. Exactly the
same proof holds for any field $K$ of characteristic different from 2.

\begin{proposition}\label{identities we need}
{\rm (i)} The factor algebra $K\langle X\rangle/T$ satisfies the identities
\[
[x_1,x_2]x_3=x_3[x_1,x_2],
\]
\[
[x_1,x_2][x_1,x_3]=0,\quad [x_1,x_2]x_4[x_1,x_3]=0,
\]
\[
[x_1,x_2][x_3,x_4]+[x_1,x_3][x_2,x_4]=0,
\]
\[
[x_1,x_2]x_5[x_3,x_4]+[x_1,x_3]x_5[x_2,x_4]=0.
\]

{\rm (ii)} The products
\begin{equation}\label{basis modulo the T-ideal}
x_{i_1}\cdots x_{i_l}[x_{j_1},x_{k_1}]\cdots[x_{j_m},x_{k_m}],
\end{equation}
$i_1\leq\cdots\leq i_l$, $k_1<j_1<\cdots<k_m<j_m$,
form a $K$-basis of $K\langle X\rangle/T$.
\end{proposition}

\begin{theorem}\label{Groebner basis of T}
Let $\text{\rm char}(K)\not=2$.
The polynomials
\[
f'_{ij}=[[x_i,x_j],x_j],\quad f''_{ij}=[x_i,[x_i,x_j]],\quad i>j,
\]
\[
g'_{ijk}=[x_i,[x_j,x_k]],\quad g''_{ikj}=[[x_i,x_k],x_j],\quad i>j>k,
\]
from (\ref{commutators of length 3 in two variables}) and
(\ref{commutators of length 3 in three variables}) and the polynomials
\begin{equation}\label{product of two commutators in two variables}
t_{ij}=[x_i,x_j][x_i,x_j],\quad i>j,
\end{equation}
\begin{equation}\label{product of two commutators in three variables}
u'_{ijk}=[x_i,x_j][x_i,x_k],\quad u''_{ijk}=[x_i,x_k][x_i,x_j],\quad i>j>k,
\end{equation}
\begin{equation}\label{product of two commutators in three variables plus A}
v'_{ijka}=[x_j,x_k]x_j^{a_j}\cdots x_{i-1}^{a_{i-1}}[x_i,x_k],
\end{equation}
\begin{equation}\label{product of two commutators in three variables plus B}
v''_{ijka}=[x_j,x_k]x_j^{a_j}\cdots x_{i-1}^{a_{i-1}}[x_i,x_j],
\end{equation}
where $i>j>k$, $a_j,\ldots,a_{i-1}\geq 0$,
\begin{equation}\label{product of two commutators in four variables plus A}
w'_{ijkla}=[x_j,x_k]x_j^{a_j}\cdots x_{i-1}^{a_{i-1}}[x_i,x_l]
+[x_j,x_l]x_j^{a_j}\cdots x_{i-1}^{a_{i-1}}[x_i,x_k],
\end{equation}
\begin{equation}\label{product of two commutators in four variables plus B}
w''_{ijkla}=[x_j,x_l]x_j^{a_j}\cdots x_{i-1}^{a_{i-1}}[x_i,x_k]
+[x_k,x_l]x_j^{a_j}\cdots x_{i-1}^{a_{i-1}}[x_i,x_j],
\end{equation}
where $i>j>k>l$, $a_j,\ldots,a_{i-1}\geq 0$,
form a minimal Gr\"obner basis with respect to the deg-lex ordering
of the T-ideal of $K\langle X\rangle$
generated by $[x_1,x_2,x_3]$.
\end{theorem}

\begin{proof}
By Proposition \ref{identities we need} (i), the polynomials
(\ref{commutators of length 3 in two variables}),
(\ref{commutators of length 3 in three variables}),
(\ref{product of two commutators in two variables}),
(\ref{product of two commutators in three variables}),
(\ref{product of two commutators in three variables plus A}),
(\ref{product of two commutators in three variables plus B}),
(\ref{product of two commutators in four variables plus A}),
(\ref{product of two commutators in four variables plus B})
belong to the T-ideal $T$ generated by $[x_1,x_2,x_3]$. Their leading terms
are obtained by deleting the commutators in the corresponding elements and are,
respectively,
\[
x_ix_jx_j,\quad x_ix_ix_j,\quad i>j,
\]
\[
x_ix_jx_k,\quad x_ix_k,x_j,\quad i>j>k,
\]
\[
x_ix_jx_ix_j,\quad i>j,
\]
\[
x_ix_jx_ix_k,\quad x_ix_kx_ix_j,\quad i>j>k,
\]
\[
x_jx_kx_j^{a_j}\cdots x_{i-1}^{a_{i-1}}x_ix_k,
\quad
x_jx_kx_j^{a_j}\cdots x_{i-1}^{a_{i-1}}x_i,x_j,\quad i>j>k,
\]
\[
x_jx_kx_j^{a_j}\cdots x_{i-1}^{a_{i-1}}x_ix_l,\quad
x_jx_lx_j^{a_j}\cdots x_{i-1}^{a_{i-1}}x_ix_k,\quad i>j>k>l,
\]
and $a_j,\ldots,a_{i-1}\geq 0$. It is easy to see that these words are pairwise different.
Clearly, the polynomial
$u'_{ijk}=[x_i,x_j][x_i,x_k]$ from
(\ref{product of two commutators in three variables})
has the same leading term as $h_{ijk} = [[x_i,x_j],[x_i,x_k]]$
from (\ref{commutators of length 4}).
Hence the set of normal words with respect to
$f'_{ij}$, $f''_{ij}$, $g'_{ijk}$, $g''_{ikj}$, $u'_{ijk}$
is the same as the one in Lemma \ref{normal words for gamma3} and we may assume that
these normal words are in the form (\ref{easy form for normal words}).
Now we want to remove the words in
(\ref{easy form for normal words}) which contain as a subword a leading word
of some $t_{ij}$, $u''_{ijk}$, $v'_{ijka}$, $v''_{ijka}$, $w'_{ijkla}$, $w''_{ijkla}$.
If $b_{ij}\geq 2$ for some $i,j$, then we remove the word using
$t_{ij}$. Hence we may assume that $b_{ij}\leq 1$. If
$b_{ik}=b_{ij}=1$ for some $i>j>k$, and $b_{i,k+1}=\cdots=b_{i,j-1}=0$, then
we use $u''_{ijk}$. Therefore, the words left in (\ref{easy form for normal words}) are
\begin{equation}\label{the words left after first reduction}
x_1^{a_1}(x_2x_1)^{\varepsilon_2}x_2^{a_2}(x_3x_{k_3})^{\varepsilon_3}x_3^{a_3}\cdots
(x_mx_{k_m})^{\varepsilon_m}x_m^{a_m},
\end{equation}
where $a_i\geq 0$, $i>k_i$, $\varepsilon_i=0,1$. Let us consider
two consecutive nonzero $\varepsilon_c$ and $\varepsilon_d$.
The corresponding monomial contains a subword
\begin{equation}\label{subwords left}
x_cx_px_c^{a_c}\cdots x_{d-1}^{a_{d-1}}x_dx_q,\quad d>c>p,d>q.
\end{equation}
If $p=q$ or $c=q$, then we use, respectively, $v'_{dcpa}$ and $v''_{dcp}$.
If $c,d,p,q$ are pairwise different, then we have the three possibilities
$p>q$, $c>q>p$, and $q>c$. The first two possibilities are excluded,
respectively, using $w'_{dcpqa}$ and $w''_{dcqpa}$.
In this way, the only subwords (\ref{subwords left}) left are for $d>q>c>p$.
Hence, we reduce the set of normal words from
(\ref{the words left after first reduction}) to the words with the condition
that for the nonzero $\varepsilon_{j_1},\ldots,\varepsilon_{j_r}$ we have
\[
k_{j_1}<j_1<k_{j_2}<j_2<\cdots<k_{j_r}<j_r.
\]
Using the correspondence $\vartheta$ from Proposition \ref{correspondence}, we obtain
that these words are in bijection with the basis elements
(\ref{basis modulo the T-ideal}) of $K\langle X\rangle/T$ which preseves the multigrading.
This implies that the polynomials
$f'_{ij}$, $f''_{ij}$, $g'_{ijk}$, $g''_{ikj}$,
$t_{ij}$, $u'_{ijk}$, $u''_{ijk}$, $v'_{ijka}$,
$v''_{ijka}$, $w'_{ijkla}$, and $w''_{ijkla}$
do form a minimal Gr\"obner basis of the T-ideal.
\end{proof}

We can state Theorem \ref{Groebner basis of T} in the following way.
We require $\vert X\vert\geq 5$ for simplification of the statement only.

\begin{theorem}
Let $\text{\rm char}(K)\not=2$, $\vert X\vert\geq 5$, and let $S$ be the semigroup
of $\text{\rm End}(K\langle X\rangle)$ generated by all endomorphisms sending $x_1,x_2,x_3,x_4$
to arbitrary elements of $X$ (allowing repetitions) and $x_5$ to products
of the form $x_1^{a_1}\cdots x_m^{a_m}$, $a_i\geq 0$.
The polynomials
\[
[[x_1,x_2],x_3],\quad
[x_1,x_2]x_5[x_3,x_4]+[x_1,x_3]x_5[x_2,x_4]
\]
form a (nonminimal) Gr\"obner $S$-basis with respect to the deg-lex ordering
of the T-ideal of $K\langle X\rangle$
generated by $[x_1,x_2,x_3]$.
\end{theorem}

\begin{remark}
(i) As in the previous section, the condition that $X$ is countable can
be replaced by the requirement that $X$ is any infinite well-ordered set.

(ii) In \cite{BML} Bokut and Makar-Limanov include in the list of the Gr\"obner basis
of the T-ideal of $K\langle x_1,x_2\rangle$ generated by $[x_1,x_2,x_3]$ the element
$(x_2x_1)^2-(x_1x_2)^2$. The evaluation of this polynomial on the Grassmann algebra
$x_1\to 1+e_1$, $x_2\to 1+e_2$ shows that $(x_2x_1)^2-(x_1x_2)^2$ does not belong to the T-ideal.
The correct Gr\"obner basis consists of the three polynomials
\[
[[x_2,x_1],x_1],\quad [x_2,[x_2,x_1]],\quad [x_2,x_1][x_2,x_1].
\]
\end{remark}

\section*{Acknowledgements}
The authors are very grateful to the anonymous referees for the
numerous suggestions which led to the improvement of the exposition and
the extension of the list of references.

\end{document}